\documentclass[11pt,a4paper, notitlepage]{article}
\usepackage{import} 
\usepackage{fullpage} 
\usepackage[utf8]{inputenc}
\usepackage[english]{babel}

\usepackage[T1]{fontenc} 
\usepackage{amssymb,amsfonts}
\usepackage{amsthm}
\usepackage{tikz-cd}
\usepackage[colorinlistoftodos]{todonotes}
\usepackage[framemethod=TikZ]{mdframed} 
\usepackage{calligra,mathrsfs}
\usepackage{amsmath}
\usepackage{mathtools}
\usetikzlibrary{shapes.geometric, arrows}
\usepackage[backend=biber,style=alphabetic,sorting=nyt, maxnames=10, maxalphanames=10]{biblatex}
\usepackage{csquotes}
\usepackage{colonequals}
\usepackage{adjustbox} 
\usepackage[colorlinks, final]{hyperref} 
\usepackage{cleveref}  
\hypersetup{
  colorlinks = true,
  breaklinks = true,
  linkcolor  = green, 
  citecolor  = blue, 
  urlcolor   = purple, 
  linktocpage = true,
}     

\newcommand{\citeafn}[1]{\citeauthorfullname{#1}}
\newrobustcmd*{\citeauthorfullname}{\AtNextCite{\DeclareNameAlias{labelname}{given-family}}\citeauthor}

\newtheorem{thm}{Theorem}[section]
\newtheorem{cor}[thm]{Corollary}
\newtheorem{prop}[thm]{Proposition}
\newtheorem{lem}[thm]{Lemma}

\theoremstyle{definition}
\newtheorem{defn}[thm]{Definition}

\theoremstyle{remark}

\newtheorem{rem}[thm]{Remark}

\numberwithin{equation}{subsection}

\newenvironment{teqn}{\begin{equation}\begin{tikzcd}\displaystyle}{\end{tikzcd}\end{equation}}
\newenvironment{teqn*}{\begin{equation*}\begin{tikzcd}\displaystyle}{\end{tikzcd}\end{equation*}}

\newenvironment{tpic*}{\begin{equation*}\begin{tikzpicture}}{\end{tikzpicture}\end{equation*}}

\DeclareMathAlphabet\matheu{U}{eus}{m}{n}
\DeclareMathAlphabet{\mathmm}{U}{mmambb}{m}{n}
\newcommand{\bb}[1]{\mathbb{#1}}
\newcommand{\fr}[1]{\mathfrak{#1}}
\newcommand{\ca}[1]{\mathcal{#1}}
\newcommand{\scr}[1]{\mathscr{#1}}
\newcommand{\eu}[1]{\matheu{#1}}

\newcommand{\bn}{\begin{enumerate}}
\newcommand{\bt}{\begin{thm}}
\newcommand{\bl}{\begin{lem}}
\newcommand{\bp}{\begin{prop}}
\newcommand{\bc}{\begin{cor}}
\newcommand{\bd}{\begin{teqn}}
\newcommand{\bud}{\begin{teqn*}}	
\newcommand{\bs}{\begin{proof}}
\newcommand{\br}{\begin{rem}}
\newcommand{\bdf}{\begin{defn}}

\newcommand{\en}{\end{enumerate}}
\newcommand{\et}{\end{thm}}
\newcommand{\el}{\end{lem}}
\newcommand{\ep}{\end{prop}}
\newcommand{\ec}{\end{cor}}
\newcommand{\ed}{\end{teqn}}
\newcommand{\eud}{\end{teqn*}}
\newcommand{\es}{\end{proof}}
\newcommand{\er}{\end{rem}}
\newcommand{\edf}{\end{defn}}
\DeclareMathOperator{\etale}{\acute{e}t}

\DeclareMathOperator{\op}{op}

\DeclareMathOperator{\spec}{Spec}
\DeclareMathOperator{\cat}{Cat}
\DeclareMathOperator{\sch}{Sch}

\DeclareMathOperator{\spa}{Spa}

\DeclareMathOperator{\fet}{\mathfrak{fet}}

\newcommand{\R}{\mathrm{R}}
\DeclareMathOperator{\sk}{sk}
\DeclareMathOperator{\cosk}{cosk}
\newcommand{\D}[1]{\mathscr{D}(#1)}

\newcommand{\tintir}{\arrow[r,shift left=0.2cm]\arrow[r]\arrow[r, shift right=0.2cm]}
\newcommand{\dutir}{\arrow[r, shift left=0.1cm]\arrow[r, shift right=0.1cm]}
\newcommand{\tir}{\arrow[r]}
\renewcommand{\i}[1]{\textit{#1}}
\newcommand{\iso}{\xrightarrow{\sim}}
\newcommand{\isom}{\arrow[r, "\sim"]}

\newcommand{\tilt}[1]{{#1^{\flat}}}
\newcommand{\untilt}[1]{{#1^{\sharp}}}


\bibliography{bibli.bib} 
\makeatletter
\newcommand*{\doublerightarrow}[2]{\mathrel{
  \settowidth{\@tempdima}{$\scriptstyle#1$}
  \settowidth{\@tempdimb}{$\scriptstyle#2$}
  \ifdim\@tempdimb>\@tempdima \@tempdima=\@tempdimb\fi
  \mathop{\vcenter{
    \offinterlineskip\ialign{\hbox to\dimexpr\@tempdima+1em{##}\cr
    \rightarrowfill\cr\noalign{\kern.5ex}
    \rightarrowfill\cr}}}\limits^{\!#1}_{\!#2}}}
\newcommand*{\triplerightarrow}[1]{\mathrel{
  \settowidth{\@tempdima}{$\scriptstyle#1$}
  \mathop{\vcenter{
    \offinterlineskip\ialign{\hbox to\dimexpr\@tempdima+1em{##}\cr
    \rightarrowfill\cr\noalign{\kern.5ex}
    \rightarrowfill\cr\noalign{\kern.5ex}
    \rightarrowfill\cr}}}\limits^{\!#1}}}
\makeatother                                                                                                                                            
\DeclareMathOperator{\idem}{Idem}
\newcommand{\etcoh}[1]{\R{\Gamma}_{\etale}(#1, \eu{G})}
\DeclareMathOperator{\spc}{\cal{S}}
\newcommand{\ccite}[2]{\cite[#2]{#1}}
\newcommand{\hjbrl}[1]{\tilde{#1}}
\newcommand{\ma}[1]{#1}
\newcommand{\inv}[1]{[\frac{1}{#1}]}
\newcommand{\aic}{valuation rings of rank $\le 1$ with algebraically closed fraction fields}

\newcommand{\aicsa}[1]{valuation rings over $#1$ of rank $\le 1$ with algebraically closed fraction fields}
\newcommand{\letperf}{$A$ be a perfectoid ring with an element $\varpi\in A$ such that $\varpi^p\mid p$ and that $A$ is $\varpi$-adically complete}
\newcommand{\lettilt}{$\tilt{A}$ be its tilt and $\tilt{\varpi}\in\tilt{A}$ be such that $\varpi^{\flat\sharp}$ is a unit multiple of $\varpi$ in $A$}
\newcommand{\giveperf}{a perfectoid ring $A$ and an element $\varpi\in A$ such that $\varpi^p\mid p$ and that $A$ is $\varpi$-adically complete}
\newcommand{\givetilt}{its tilt $\tilt{A}$ with an element $\tilt{\varpi}\in\tilt{A}$ such that $\varpi^{\flat\sharp}$ is a unit multiple of $\varpi$ in $A$}
\newcommand{\bcpair}{$U\subset\spec A$ (resp.~$\tilt{U}\subset\spec\tilt{A}$) is an open containing $\spec(A[\frac{1}{\varpi}])$, (resp.~containing $\spec(\tilt{A}[\frac{1}{\tilt{\varpi}}])$,) such that the closed subsets $Z\colonequals \spec A\setminus U$ and $\tilt{Z}\colonequals\spec\tilt{A}\setminus\tilt{U}$ agree under the isomorphism induced by \eqref{diag:identification}}

\title{\uppercase{\textbf{\large{Tilting Correspondences of Perfectoid Rings}}}}
\author{\uppercase{Arnab Kundu}\let\thefootnote\relax\footnote{Universit\'e Paris-Saclay, Laboratoire de Math\'ematiques d’Orsay, F-91405, Orsay, France\newline\hspace*{0.48cm} Email: arnab.kundu@universite-paris-saclay.fr.}}
\date{}
\usepackage{doi}
\usepackage{microtype} 
\begin{document}
\maketitle
\abstract{In this article, we present an alternate proof of a vanishing result of étale cohomology on perfectoid rings due to Česnavičius and more recently proved by a different approach by Bhatt and Scholze. To establish that, we prove a tilting equivalence of étale cohomology of perfectoid rings taking values in commutative, finite étale group schemes. On the way, we algebraically establish an analogue of the tilting correspondences of Scholze, between the category of finite étale schemes over a perfectoid ring and that over its tilt, without using tools from almost ring theory or adic spaces.}
\tableofcontents
\let\thefootnote\relax\footnote{Date: October 2020}
\section{Introduction}
Our goal in this article is to simplify the proof of the following result of \citeauthor{ces-brauer}. This theorem plays a central role in the proof of the purity of the Brauer group, notably by demonstrating that the $p$-primary torsion part of the Brauer group $\textstyle H^2_{\etale}(A[\frac{1}{p}],\bb{G}_m)[p^\infty]$ for a perfectoid ring $A$ vanishes. 
\bt[\ccite{ces-brauer}{Thm.~4.10}, cf.~\ccite{bhatt-scholze}{Thm.~11.1}, see \Cref{cor:main}]\label{thm:ces-brauer}
Let $p$ be a prime and let $A$ be a $\bb{Z}_p$-algebra such that it is a perfectoid ring with an element $\varpi\in A$ such that $\varpi^p\mid p$ and that $A$ is $\varpi$-adically complete. Then, for a commutative, finite étale $A[\frac{1}{\varpi}]$-group scheme $G$ of $p$-power order, we have, for all $i\ge 2$, \bud\textstyle H^i_{\etale}(A[\frac{1}{\varpi}],G)=0.\eud
\et
We follow the definition of a perfectoid ring from \cite{bms1} (see \Cref{defn:integral-perfd}). We remark that $H^1_{\etale}(A[\frac{1}{\varpi}], G)$ can be computed using the prismatic Dieudonn\'e module of $G$ (see \ccite{ces-scholze}{Thm.~4.1.8}). 
The statement of \Cref{thm:ces-brauer} is a mild generalisation of \cite[Thm.~4.10]{ces-brauer}. Indeed, the beginning of \textsection\ref{section:arc-descent} implies that we can find $\pi\in A$ such that $A$ is $\pi$-adically complete and that $\pi^p=p$; consequently, $A[\frac{1}{\pi}]=A[\frac{1}{p}]$ and for all $i\ge 2$, \bud\textstyle H^i_{\etale}(A[\frac{1}{p}],G)=0.\eud
This slight improvement is possible since our results do not depend on the almost purity theorem of \cite{kedlaya-liu}. The proof given in \cite{ces-brauer} uses a non-noetherian version of a result of \citeauthor{huber96} from \cite{huber96} to compare the \'etale cohomology of $\textstyle A[\frac{1}{p}]$ with coefficients in $G$ with the \'etale cohomology of the associated adic space $\textstyle\spa(A[\frac{1}{p}], A)$ with coefficients in the group associated to $G$, thus transferring the problem to the world of perfectoid spaces studied by \citeauthor{scholze-thesis} in \cite{scholze-thesis} and \citeauthor{kedlaya-liu} in \cite{kedlaya-liu}. 
 In this article, we propose a proof of \Cref{thm:ces-brauer} by using a similar strategy inspired from the work of \citeauthor{ces-scholze} in \cite{ces-scholze}, by replacing the almost purity theorem with `algebraic' tilting results from op.~cit. By bootstrapping from the arguments in \cite[Thm.~2.2.7]{ces-scholze}, we deduce a non-constant coefficient version of the tilting result of \'etale cohomology of perfectoid rings (see \Cref{thm:B}). The details of the proof of \Cref{thm:ces-brauer} using this non-constant coefficient version of the tilting result are given in \Cref{cor:main}.
\par \Cref{thm:ces-brauer} is obtained as a consequence of the isomorphism (\Cref{thm:B}) between the \'etale cohomology group appearing in \Cref{thm:ces-brauer} and the \'etale cohomology group of the corresponding perfect $\bb{F}_p$-algebra `tilt' taking values in a commutative, finite étale group scheme $\tilt{G}$ of $p$-power order. Indeed, we may apply the vanishing of the \'etale cohomology of the $\bb{F}_p$-algebra. More precisely, by \cite[Ex.~X Thm.~5.1]{sga4iii}, the $p$-cohomological dimension of an affine noetherian $\bb{F}_p$-algebra is $\le 1$, and consequently, by limit arguments, the \'etale cohomology of degree $i\ge 2$ of any commutative, finite étale group scheme of $p$-power order over an $\bb{F}_p$-algebra vanishes. 
The existence of the `tilt' $\tilt{G}$ of $G$ will be shown by applying \Cref{thm:A}, which can be seen as an algebraic analogue of tilting results on perfectoid Banach $K$-algebras over a perfectoid field $K$ as in \cite{scholze-thesis} or that on perfectoid Banach $\bb{Q}_p$-algebras as in \cite{kedlaya-liu}. 
\bt\label{thm:A}
Let $A/\bb{Z}_p$ be a perfectoid ring with an element $\varpi\in A$ such that $\varpi^p\mid p$ and that $A$ is $\varpi$-adically complete, and let $\tilt{A}$ be its tilt and $\tilt{\varpi}\in\tilt{A}$ be such that $\varpi^{\flat\sharp}$ is a unit multiple of $\varpi$ in $A$. 
Then, there is an equivalence, functorial in $A$, between the categories of finite étale algebras \bd\tag{$\ast$}\label{diag:thmA}\textstyle\fet/A[\frac{1}{\varpi}]\cong\fet/\tilt{A}[\frac{1}{\tilt{\varpi}}].\ed 
\et
For a ring $R$, the notion of a (co-commutative) Hopf $R$-algebra is dual to the notion of a (commutative) affine $R$-group scheme, i.e., they are $R$-algebras so that their spectra have the structure of (commutative) affine $R$-group schemes. 
Starting with a commutative, finite \'etale group $A[\frac{1}{\varpi}]$-scheme $G$, we observe that, by \eqref{diag:thmA}, its coordinate ring $\ca{O}(G)$, which is a co-commutative, finite \'etale Hopf $A[\frac{1}{\varpi}]$-algebra, has a tilt $\tilt{\ca{O}(G)}$, which a priori is a finite \'etale $\tilt{A}[\frac{1}{\tilt{\varpi}}]$-algebra. However, the structure of a (co-commutative) Hopf algebra, expressible in terms of diagrams involving $\textstyle A[\frac{1}{\varpi}]$, $\ca{O}(G)$ and $\displaystyle\ca{O}(G)\otimes_A\ca{O}(G)$, gets transferred by \eqref{diag:thmA} to $\tilt{\ca{O}(G)}$, whence we get the the required commutative, finite étale $\tilt{A}[\frac{1}{\tilt{\varpi}}]$-group scheme $\tilt{G}\coloneqq\spec\tilt{\ca{O}(G)}$, as shown in the following theorem.
\bt[see \Cref{thm:tilting-U-fet}]\label{thm:B}
We assume notations of \Cref{thm:A}. For commutative group schemes $\textstyle G\in\fet/A[\frac{1}{\varpi}]$ and $\textstyle\tilt{G}\in\fet/\tilt{A}[\frac{1}{\tilt{\varpi}}]$ that are identified under the tilting correspondence \eqref{diag:thmA}, 
we have an identification, functorial in $A$ and $G$, \bd\tag{$\ast\ast$}\label{diag:thmB}\textstyle R\Gamma_{\etale}(A[\frac{1}{\varpi}], G)\cong R\Gamma_{\etale}(\tilt{A}[\frac{1}{\tilt{\varpi}}], \tilt{G}).\ed
\et
We have the following corollary to \Cref{thm:A}.
\bc\label{cor:thmA}
Let $R$ be a perfectoid Banach $K$-algebra as in \cite[Defn.~5.1]{scholze-thesis}, where $K$ is a perfectoid field (resp.~a perfectoid Banach $\bb{Q}_p$-algebra as in \cite[Defn.~3.6.1]{kedlaya-liu}). Choose an element $\varpi\in R^\circ$ such that $\varpi^p\mid p$ and that $R^\circ$ is $\varpi$-adically complete, and choose $\tilt{\varpi}\in\tilt{{R^\circ}}$ such that $(\varpi^{\flat\sharp})=(\varpi)$. Then, there is a functorial in $R$ equivalence of categories, \bud\textstyle\fet/R^\circ[\frac{1}{\varpi}]\cong\fet/R^{\flat\circ}[\frac{1}{\tilt{\varpi}}].\eud 
\ec
\par A key ingredient in the proof of \Cref{thm:B} is the theory of the arc topology (\Cref{defn:arc-topology}) developed by \citeauthor{bhatt-mathew-arc} in \cite{bhatt-mathew-arc}. As shown by them, the cohomology of an \'etale sheaf of torsion abelian groups satisfies (hyper)descent in the arc topology (see \Cref{defn:hyperdescent} for the definition of hyperdescent and \Cref{prop:fet-arc-descent} for a precise statement). 
Moreover, they proved that, given an \'etale sheaf $\eu{G}$ of torsion abelian groups, the functor $R'\mapsto\R{\Gamma}_{\etale}(\hat{R'}[\frac{1}{\varpi}], \eu{G})$ satisfies $\varpi$-complete arc hyperdescent (the $\varpi$-complete arc topology defined in \Cref{defn:w-arc-top} is a slight refinement of the arc topology and is better suited to the study of perfectoid rings). 
This result (\Cref{rem:arc-w-descent}) and the fact that any perfectoid ring has a $\varpi$-complete arc hypercover given by a special class of perfectoid rings (see \Cref{lem:arc-hypercover2}) reduces the proof to showing the equivalence for this special class of perfectoid rings, where it is possible to give a direct proof. It must be remarked, however, that except for the aforementioned class of special perfectoid rings, we do not know any direct morphism that establishes an identification \eqref{diag:thmB}. This failure is related to the fact that the `tilting functor' (cf.~\Cref{prop:alg-tilting}), which localises in the analytic topology in the adic case, does not localise in the Zariski topology in the algebraic case. 
\subsection*{Notations and Conventions}
We fix a prime integer $p$. The term \textit{rank} shall denote the Krull dimension of a valuation ring. Given a ring $A$ with an element $\varpi\in A$ and an $n\ge 1$, we define $A\langle\varpi^n\rangle$ to be the kernel of multiplication by $\varpi^n$ map in $A$; these kernels form an increasing system with union $A\langle\varpi^{\infty}\rangle$. Given a ring $R$, the category of \textit{schemes over $R$} will be denoted by $\sch_R$ and its subcategory of \textit{quasi-compact and quasi-separated $R$-schemes} will be denoted by $\sch^{\mathrm{qcqs}}_R$. Given a simplicial object $X_{\bullet}$ we shall denote the \textit{$n$-th component} by $X_n$ and the \textit{$n$-truncation} of the object by $X_{\le n}$. For any category $\scr{C}$ and an object $X\in\scr{C}$, the \textit{slice category over $X$} will be denoted by $\scr{C}/X$. The $2$-category of small $1$-categories will be denoted $\cat$.
 
\subsection*{Acknowledgement}
This paper is a product of the author's Master's Thesis \cite{master-thesis} done under the guidance of ~\citeafn{ces-brauer}. His role has been crucial as his work served as the inspiration for this project. The author expresses gratitude for his constant encouragement, unwavering support and careful readings which has improved this manuscript in several aspects. The author would like also to thank Thibault Alexandre, Vincent Bouis, Arnaud Eteve, Ning Guo, Kazuhiro Ito, Hiroki Kato, Ariane Mézard, Mohamed Moakher, Shubhodip Mondal, Shravan Patankar, Maxime Ramzi, Benoît Stroh and Yifei Zhao for stimulating conversations and comments, and Marc Hoyois for suggesting the paper \cite{ehik}. This project received partial funding from the Fondation Sciences Mathématiques de Paris' PGSM Master's scholarship and the European Research Council under the European Union’s Horizon 2020 research and innovation program (grant agreement No. 851146).
\section{\texorpdfstring{$\varpi$}{w}-complete Arc (hyper)sheaves \texorpdfstring{$\fet$}{fet} and \texorpdfstring{$\R\Gamma_{\etale}(-,\eu{G})$}{RGamma(-,G)}}\label{section:arc-descent}
Our goal of this section is to deduce \Cref{rem:arc-w-descent} from \Cref{prop:fet-arc-descent} (following \cite{bhatt-mathew-arc}). 
It is essential in this article to use $\infty$-categorical tools, and therefore, we refer the reader to \cite{lurie-htt}. 
\begin{defn}[\ccite{stacks-project}{\href{https://stacks.math.columbia.edu/tag/049J}{Tag 049J}}]\label{defn:lfp}
For a ring $R$ and an $\infty$-category $\scr{C}$ with all colimits, we say that a functor $\eu{F}\colon\sch^{\op}_R\to\scr{C}$ is \textit{locally of finite presentation} (`finitary' in \cite{bhatt-mathew-arc}), if whenever $\{S_{\alpha}, f_{\alpha\beta}\}_{\alpha,\beta\in I}$ is an inverse system of quasi-compact and quasi-separated $R$-schemes indexed by a cofiltered partially ordered set $I$, such that the transition map $f_{\alpha\beta}\colon S_{\alpha}\to S_{\beta}$ is affine for each $\alpha,\beta\in I$, the map $\textstyle\mathrm{colim}_{\alpha}\left(\eu{F}(S_\alpha)\right)\xrightarrow{\sim}\eu{F}(\lim_{\alpha} S_{\alpha})
$ is an equivalence in $\scr{C}$. 
\end{defn}
\begin{defn}[\ccite{stacks-project}{\href{https://stacks.math.columbia.edu/tag/01G5}{Tag 01G5}}, \ccite{lurie-htt}{Defn.~6.5.3.2}]\label{defn:hypercover}
Given a category $\scr{X}$ with finite limits, an element $X\in\scr{X}$, and a Grothendieck topology $\tau$ on $\scr{X}$, a simplicial object $X_{\bullet}$ in $\scr{X}/X$ is a \textit{hypercover} of $X$ if 
\bn
\item the morphism $X_0\to X$ is a $\tau$-covering, and
\item for every $n\ge 1$, the morphism $X_{n}\to(\cosk_n(\sk_{n-1}(X_{\bullet})))_n$ is a $\tau$-cover. 
\en 
\end{defn}
For example, given a $\tau$-cover $Y\to X$, we have the \v{C}ech nerve hypercover 
\bud(\cdots\tintir &Y\times_X Y\dutir & Y)\text{ of }X.\eud


\begin{defn}[cf.~\ccite{lurie-htt}{\textsection6.5.2}]\label{defn:hyperdescent}
Given a ring $R$, a Grothendieck topology $\tau$ on $\sch_R$, and an $\infty$-category $\scr{C}$ with all limits, a functor $\eu{F}\colon\sch^{\op}_R\to\scr{C}$ is said to be a $\tau$\textit{-sheaf} or to satisfy \textit{$\tau$-descent}, if $\eu{F}$ carries finite coproducts of schemes to products in $\scr{C}$ and for every $\tau$-cover $Y\to X$, the map \bud \eu{F}(X)\isom & \lim(\eu{F}(Y)\dutir & \eu{F}(Y\times_X Y)\tintir & \cdots)\text{ is an equivalence.} \eud A $\tau$-sheaf $\eu{F}$ is said to be a \i{$\tau$-hypersheaf} or to satisfy \i{$\tau$-hyperdescent}, if for every hypercover 
$X_{\bullet}\coloneqq(\cdots\triplerightarrow{}X_1\doublerightarrow{}{}X_0)$ of $X$ in the $\tau$-topology, the map \bud
\eu{F}(X) \isom & \lim(\eu{F}(X_0) \dutir & \eu{F}(X_1) \tintir & \cdots)\text{ is an equivalence.}\eud
\end{defn} 
Given an $\infty$-category $\scr{C}$ and two objects $C,D\in\scr{C}$, \cite[Defn.~1.2.2.1]{lurie-htt} associates a Kan complex $\mathrm{Map}_{\scr{C}}(C,D)$, called the `mapping space' from $C$ to $D$. This is a generalisation of the set of morphisms from $C$ to $D$ in an $1$-category, obtained by admitting homotopies between morphisms. For $n\ge -1$, an object $C\in\mathscr{C}$ is said to be \textit{$n$-truncated} if $\mathrm{Map}_{\scr{C}}(D,C)$ is $n$-truncated, for all objects $D\in\scr{C}$ (see \cite[Defn.~5.5.6.1]{lurie-htt}), i.e., if for all $k>n$ the homotopy groups $\pi_k(\mathrm{Map}_{\scr{C}}(D,C))$ vanish; dually $C$ is \textit{$n$-cotruncated} (resp., cotruncated) if it is $n$-truncated (resp., $m$-cotruncated for some $m\ge -1$) in $\scr{C}^{\op}$. A cocomplete $\infty$-category $\scr{C}$ is said to be \textit{generated under colimits by cotruncated objects} if the class of cotruncated objects forms a set, and any object $C\in\scr{C}$ can be written as a colimit of cotruncated objects (cf.~\cite[Defn.~3.1.4]{ehik}). Examples include `nice' cocomplete $n$-categories (in which every object is $n$-truncated), where `nice' refers to some set theoretic finiteness condition to make the category small enough, and the derived $\infty$-category $\D{\bb{Z}}^{\ge 0}$ of bounded below by $0$ complexes of abelian groups (in which the perfect complexes form the set cotruncated objects by \cite[Warning 1.2.1.9]{lurie-algebra} and \cite[\href{https://stacks.math.columbia.edu/tag/07VI}{Tag 07VI}]{stacks-project}; also cf.~\cite[Eg.~3.5(1)]{bhatt-mathew-arc}). 
\bp[\ccite{ehik}{Lem.~3.1.7}]\label{prop:descent->hyperdescent}
For a ring $R$, a Grothendieck topology $\tau$ on $\sch_R$, and an $\infty$-category $\scr{C}$ generated under colimits by cotruncated objects, any sheaf $\eu{F}\colon\sch_R^{\op}\to\scr{C}$ is automatically a hypersheaf.
\ep
Following notations of op.~cit., we denote the $\infty$-category of Kan complexes as $\spc$ and for each $n\ge -1$, its full $\infty$-subcategory of $n$-truncated Kan complexes as $\spc_{\le n}$. In loc.~cit.~the source of $\eu{F}$ is assumed to be an `$\infty$-topos', however, the same proof works in the case of an ordinary Grothendieck site, which is sufficient for our purposes. We now present a sketch of the proof of \Cref{prop:descent->hyperdescent} following arguments from op.~cit.~and \cite{lurie-htt}.
\bs The property of a sheaf (or a hypersheaf) can be tested by applying the Yoneda emedding; whence, since $\scr{C}$ is generated under colimits by cotruncated objects, $\eu{F}$ has the property if and only if for every cotruncated object $C\in\scr{C}$, the functor $\mathrm{Map}_{\scr{C}}(C, \eu{F}(-))\colon\sch_R^{\op}\to\spc_{\le n}$ has the same property. 
The latter is a truncated object in the category of sheaves on $\sch_R$ taking values in $\spc$, and therefore, by \cite[Lem.~6.5.2.9]{lurie-htt}, it is a hypersheaf. 
\es
We have inherently used \cite[Cor.~6.5.3.13]{lurie-htt}, which shows that the notion of a hypersheaf coincides with the notion of `hypercompleteness' (see \cite[\textsection6.5.2]{lurie-htt} for the definition).

\par The following topology was defined by \citeauthor{rydh} in \cite[\textsection 1]{rydh} as `universally submersive topology' to study descent properties of étale morphisms. It was further developed by \citeauthor{bhatt-mathew-arc} in \cite{bhatt-mathew-arc}, where they proved that the étale cohomology satisfies arc descent (see \Cref{prop:fet-arc-descent}). 
\begin{defn}[\ccite{bhatt-mathew-arc}{Defn.~1.2}]\label{defn:arc-topology}
A morphism $X'\to X$ of schemes is an \textit{arc cover} if for every rank $\le 1$ valuation ring $V$ and every morphism $\spec V\to X$, there are a faithfully flat extension $V\to V'$ of valuation rings and a morphism $\spec V'\to X'$ that lifts the composition $\spec V'\to\spec V\to X$ to a commutative square
\bd\label{diag:arc-topology}
\spec V' \arrow[r] \arrow[d] & X' \arrow[d] \\
\spec V \arrow[r] & X.
\ed 
\end{defn}
If we remove the rank $\le 1$ hypothesis from \Cref{defn:arc-topology} and consider any valuation ring $V$, we obtain the \textit{$v$-topology} of \cite[Defn.~1.1]{bhatt-mathew-arc} or the `universally subtrusive topology' of \cite[Defn.~2.2]{rydh}. In the proof of the following we shall need to use a result \cite[Thm.~3.12 and Rem.~3.13]{rydh} of \citeauthor{rydh} which says that any $v$-cover of schemes has a refinement given by a composition of a Zariski open covering and a proper surjective morphism of finite presentation, and therefore, a Zariski sheaf which satisfies descent for proper surjective morphisms of finite presentation is a $v$-sheaf. 
\begin{prop}[\ccite{bhatt-mathew-arc}{Thm.~5.6(2), Thm.~5.4}]\label{prop:fet-arc-descent}
Given a ring $R$, the functor $\fet$ that associates to a quasi-compact and quasi-separated $R$-scheme $X$ the category of finite étale $X$-schemes satisfies arc hyperdescent; moreover, the functor given by $X\mapsto\R{\Gamma}_{\etale}(X,\eu{G})$ satisfies arc hyperdescent. 
\end{prop}
We need the following lemma to prove the above proposition.
\bl[\ccite{stacks-project}{\href{https://stacks.math.columbia.edu/tag/09ZL}{Tag 09ZL}}, \ccite{gabber_1994}{Thm.~1}]\label{lem:lifting-et}
Given a ring $A$ with an ideal $I\subset A$ such that $A$ is $I$-henselian, there is an equivalence $\fet/A\iso\fet/(A/I)$; moreover, given an étale sheaf $\eu{G}$ on $A$ of torsion abelian groups, the morphism $\R{\Gamma}_{\etale}(A,\eu{G})\iso\R{\Gamma}_{\etale}(A/I,\eu{G})$ is an isomorphism \el 

\bs[Proof of \Cref{prop:fet-arc-descent}]
This proof is the same as the proofs in \cite{bhatt-mathew-arc}. We shall prove that $\fet$ satisfies arc descent, since then it automatically satisfies arc hyperdescent by \Cref{prop:descent->hyperdescent}. Indeed, $\scr{C}$ can taken to be the essentially small $2$-category generated under colimits by the categories of finite \'etale $X$-schemes, where $X$ is a finite type $R$-scheme.
\par By \cite[Thm.~4.1]{bhatt-mathew-arc}, since $\fet$ is locally of finite presentation, it is enough to prove that $\fet$ is a $v$-sheaf, and that for any valuation ring with an algebraically closed fraction field and any prime ideal $\fr{p}\subset V$, the square \bd\label{diag:prop-fet-arc-descent} \fet/V\tir\arrow[d] & \fet/(V/\fr{p})\arrow[d] \\ \fet/V_{\fr{p}}\tir & \fet/\kappa(\fr{p})\ed is cartesian. Indeed, \eqref{diag:prop-fet-arc-descent} is cartesian because each ring among $V/\fr{p},V_{\fr{p}}$ and $\kappa(\fr{p})$ is a valuation ring with an algebraically closed fraction field, implying that each is a strictly henselian local ring. Since $V\to V/\fr{p}$ is a morphism between henselian local rings with the same residue field, the top horizontal morphism is an identification by \Cref{lem:lifting-et}, and similarly, the bottom horizontal morphism is an isomorphism.
\par Thus, it reduces to prove that $\fet$ is a $v$-sheaf: by \cite[Thm.~3.12 and Rem.~3.13]{rydh}, it suffices to show that it satisfies descent for proper surjective morphisms of finite presentation, which follows from \cite[Exp.~IX Thm.~4.12]{sga1}. 
\par The proof for the second functor has the same structure as the proof for the first. The proper base change theorem \cite[Exp.~XII Thm.~5.1]{sga4iii} and \ccite{bhatt-mathew-arc}{Lem.~5.1} show that the functor satisfies descent for proper surjective morphisms of finite presentation (see the proof of \cite[Thm.~5.4]{bhatt-mathew-arc} for details), and hence, by \cite[Thm.~3.12 and Rem.~3.13]{rydh}, it satisfies $v$-descent. The square analogous to \eqref{diag:prop-fet-arc-descent} is cartesian because the horizontal morphisms are isomorphism by \Cref{lem:lifting-et}; the lemma being applicable as the corresponding horizontal morphisms are between henselian rings with the same residue fields. The result \cite[\href{https://stacks.math.columbia.edu/tag/09YQ}{Tag 09YQ}]{stacks-project} proves that the functor is locally of finite presentation, consequently, the functor is an arc sheaf thanks to \cite[Thm.~4.1]{bhatt-mathew-arc}, and finally by \Cref{prop:descent->hyperdescent}, where we take $\scr{C}=\D{\bb{Z}}^{\ge 0}$, an arc hypersheaf. 
\es

\begin{defn}[\ccite{bhatt-mathew-arc}{Defn.~6.14}, \ccite{ces-scholze}{\textsection 2.2.1}]\label{defn:w-arc-top}
Given a ring $A$ with an element $\varpi\in A$, a morphism $A\to A'$ of rings is a \textit{$\varpi$-complete arc cover} if for any rank $\le 1$ and $\varpi$-adically complete valuation ring $V$ and any morphism $A\to V$, there are a faithfully flat extension $V\to V'$ of valuation rings and a morphism $\spec V' \to X'\colonequals \spec A' $ lifting the composition $\spec V' \to\spec V \to X\colonequals\spec A$ to a commutative square as in \eqref{diag:arc-topology}.
\end{defn}
\br\label{rem:redn}
For a ring $A$, an element $\varpi\in A$, and a $\varpi$-complete arc cover $A\to A'$, the reduction $A/\varpi\to A'/\varpi$ is an arc cover. Conversely, an arc cover $A\to A'$ is a $\varpi$-complete arc cover (see \cite[\href{https://stacks.math.columbia.edu/tag/090T}{Tag 090T}]{stacks-project}).
\er

\br Given a ring $A$ with an element $\varpi\in A$ such that $A$ is $\varpi$-adically complete, the functor $\fet$ on the category of $\varpi$-adically complete $A$-algebras taking such an $A$-algebra $A'$ to the category of finite \'etale $A'$-algebras is a $\varpi$-complete arc sheaf by \Cref{rem:redn} and \Cref{lem:lifting-et} (and hence a $\varpi$-complete arc hypersheaf by \Cref{prop:descent->hyperdescent}), similarly, given an \'etale sheaf $\eu{G}$ on $A$ of torsion abelian groups, the functor $A'\mapsto\R{\Gamma}_{\etale}(A',\eu{G})$ on the same category is a $\varpi$-complete arc hypersheaf. 
\er
\bp[cf.~\ccite{mathew}{Thm.~5.19}, \ccite{bhatt-mathew-arc}{Cor.~6.17}]\label{rem:arc-w-descent}
Given a ring $R$ with an element $\varpi\in R$, the functor taking an $R$-algebra $R'$, with $\varpi$-adic completion $\hat{R'}$, to the category of finite \'etale $\textstyle \hat{R'}[\frac{1}{\varpi}]$-algebras satisfies $\varpi$-complete arc hyperdescent. Moreover, given an \'etale sheaf $\eu{G}$ on $R$ of torsion abelian groups, the functor $R'\mapsto\R{\Gamma}_{\etale}(\hat{R'}[\frac{1}{\varpi}], \eu{G})$ satisfies $\varpi$-complete arc hyperdescent. 
\ep
See \ccite{mathew}{Thm.~5.19} for the case of finite \'etale algebras and \ccite{bhatt-mathew-arc}{Cor.~6.17} for the case of the \'etale cohomology. In either case, it reduces to prove arc descent for the respective functors. Notably, the case of the \'etale cohomology is easier, given that there is a cartesian square \bd
			\label{diag:rem-arc-w-descent}
				\etcoh{R}\arrow[r]\arrow[d] & \etcoh{R\inv{\varpi}}\arrow[d]
				\\ \etcoh{\hat{R}} \arrow[r] & \etcoh{\hat{R}\inv{\varpi}},
	\ed which, since $\D{\bb{Z}}^{\ge 0}$ is a stable $\infty$-category (\ccite{lurie-algebra}{Defn.~1.1.1.9}), is a co-cartesian square (see \ccite{lurie-algebra}{Prop.~1.1.3.4}). Consequently, it suffices to show that each of the other three functors appearing in \eqref{diag:rem-arc-w-descent} are arc sheaves, which follows from \Cref{prop:fet-arc-descent} and \Cref{lem:lifting-et}. But in the case of finite \'etale schemes, loc.~cit.~does not apply and thus, the above proof can not be naively adapted to work. 
\section{Tilting for \texorpdfstring{$\fet$}{fet} and \texorpdfstring{$\R{\Gamma}_{\etale}(-,\eu{G})$}{RGamma(-,G)} over Perfectoid Rings}\label{section:applications}
In the first part of \textsection\ref{section:applications} we recall the basics of tilting of perfectoid rings (following \cite{ces-scholze}) which form the base of the proof of \Cref{thm:tilting-U-fet}. As a corollary, we have the generalisation of \ccite{ces-brauer}{Thm. 4.10} stated in \Cref{cor:main}.
\par Given a $\bb{Z}_p$-algebra $A$, we define the \textit{tilt} $\tilt{A}\colonequals\lim_{x\mapsto x^p}A/(p)$.\label{lem:3.2}
Given a ring $A$ and an element $\varpi\in A$ such that $\varpi\mid p$ and $A$ is $\varpi$-adically complete, \cite[Lem.~3.2(i)]{bms1} implies that there is a multiplicative monoidal isomorphism \bd\label{diag:lem-3.2}\textstyle\lim_{x\mapsto x^p}A\xrightarrow{\sim}\tilt{A}.\ed The \textit{untilt} morphism is the projection onto the first factor $\sharp\colon\tilt{A}\to A$. It can be written as a composition of the Teichm\"{u}ller map $a\mapsto [a]$ and the canonical morphism $\theta$ from the $p$-typical Witt vectors of $\tilt{A}$ to $A$ defined by $[a]\mapsto\untilt{a}$.

\begin{defn}[\ccite{bms1}{Defn.~3.5}]\label{defn:integral-perfd}
A ring $A$ with tilt $\tilt{A}$ is called \textit{perfectoid} (also refered to as `integral perfectoid' by some authors) if there exists a $\varpi\in A$ such that $\varpi^p\mid p$ and such that $A$ is $\varpi$-adically complete and
\bn
\item the $p$-power morphism $x\mapsto x^p$ in $A/(\varpi)\to A/(\varpi^p)$ is surjective, and
\item\label{pt:2} the kernel of the map $\theta\colon W(\tilt{A})\to A$ (defined above) is principal.
\en
\end{defn}
Any perfect $\bb{F}_p$-algebra is perfectoid because we may take $\varpi=0$. 
\par Let \letperf. We recall from \cite[\textsection 2.1.3]{ces-scholze} that $A$ has bounded $\varpi$-torsion; more precisely, we have $A\langle\varpi\rangle=A\langle\varpi^{\infty}\rangle$. Additionally, by \cite[Lem.~3.9]{bms1}, there exists a unit $u\in A$ such that the element $u\varpi$ admits compatible $p$-power roots in $A$, that is, there exists an element $\tilt{\varpi}$ with untilt $u\varpi$. Then, by the proof of \cite[Lem.~3.10]{bms1}, there is an isomorphism \bd\label{diag:identification}A/(\varpi)\cong\tilt{A}/(\tilt{\varpi}).\ed Moreover, since $\tilt{A}$ is perfect, arguing as in the proof of \cite[Cor.~3.2.3]{bhatt-perfectoid}, the above isomorphism implies that $\tilt{A}$ is $\tilt{\varpi}$-adically complete. 
The isomorphism \eqref{diag:lem-3.2} of multiplicative monoids extends (thanks to \cite[\textsection 2.1.7]{ces-scholze}) to an isomorphism \bd\label{pt:3}\textstyle\tilt{A}[\frac{1}{\tilt{\varpi}}]\xrightarrow{\sim}\lim_{x\mapsto x^p}A[\frac{1}{\varpi}].\ed 
\begin{prop}[\ccite{ces-scholze}{Prop. 2.1.9}]\label{prop:alg-tilting}
Given \giveperf, and \givetilt; we have an equivalence between the categories of $\varpi$-adically complete $A$-algebras which are \aic~and the category of $\tilt{\varpi}$-adically complete $\tilt{A}$-algebras which are \aic.
\end{prop}
In fact, op.~cit.~proves a stronger result, by showing an equivalence between the categories of perfectoid rings over $A$ and the same over $\tilt{A}$ (cf.~\cite[Thm.~5.2]{scholze-thesis} and \cite[Thm.~3.6.5]{kedlaya-liu}). It is to be noted, however, for the sake of clarity, that $\varpi$-adically complete \aicsa{A} are perfectoid rings (resp., \aic~of characteristic $p$ are perfect rings) by \cite[\textsection 2.1.1]{ces-scholze}. 

\begin{lem}[\ccite{ces-scholze}{Lem.~2.2.2, Lem.~2.2.3}]\label{lem:arc-hypercover2}
Let \letperf, and let \lettilt. Then, there exists a $\varpi$-complete arc cover $A\to A'$ such that $A'=\prod_{i\in I}V_i$, where $I$ is an indexing set, and $V_i$ is a $\varpi$-adically complete valuation ring over $A$ of rank $\le 1$ with an algebraically closed fraction field for each $i\in I$, and $\tilt{A}\to\tilt{A'}$ is a $\tilt{\varpi}$-complete arc cover.
\end{lem}
\begin{lem}[\ccite{ces-scholze}{Lem.~2.2.4}]\label{lem:aic-v-fet}
Given a collection of valuation rings $\{V_i\}_{i\in I}$ with algebraically closed fraction fields, any \'etale cover of a quasi-compact open $U\subset \spec(\prod_{i\in I}V_i)$ has a section, and in particular, any finite étale $U$-scheme is a finite disjoint union of subsets $T\subset U$ which are both open and closed. 
\end{lem} 
\begin{prop}[Beauville--Laszlo]\label{prop:beauville-laszlo}
Let $R$ be a ring with an element $\varpi\in R$ such that $R$ has bounded $\varpi$-torsion (i.e., there exists an $n\ge 1$ such that $R\langle\varpi^{\infty}\rangle=R\langle\varpi^n\rangle$) and $S$ be a flat, $\varpi$-henselain $R$-algbebra with $\varpi$-adic completion $\hat{S}$. Given a quasi-compact open $\spec(R\inv{\varpi})\subset U\subset\spec(R)$ for which $\spec(R\inv{\varpi}\times S)$ is an fpqc cover, and an arc sheaf $\eu{F}$ locally of finite presentation on $R$ such that for every $R$-henselian pair $(A,I)$, the map $\eu{F}(A)\hookrightarrow\eu{F}(A/I)$ is injective, we have a cartesian diagram 
%
	\bd
			\label{diag:prop-beauville-laszlo-idem}
				\eu{F}(U)\arrow[r]\arrow[d] & \eu{F}(R\inv{\varpi})\arrow[d]
				\\ \eu{F}(\hat{S}) \arrow[r] & \eu{F}(\hat{S}\inv{\varpi}).
	\ed In particular, the above holds for the functor of idempotents and the functor of finite \'etale algebras, and given a torsion étale sheaf $\eu{G}$ on $R$, for the functor $\eu{F}=\R{\Gamma}_{\etale}(-,\eu{G})$.
\end{prop} 

\begin{proof}
By fpqc descent, we have an equivalence \begin{teqn*}\eu{F}(U)\tir &\lim(\eu{F}(R\inv{\varpi})\times \eu{F}(S)\dutir &\eu{F}(S\inv{\varpi})).\end{teqn*} Thanks to \cite[Thm.~2.1.16]{bouthier-ces}, the hypothesis implies that $\eu{F}(S\inv{\varpi})\hookrightarrow\eu{F}(\hat{S}\inv{\varpi})$. It remains to show that there is an equivalence \begin{teqn*}\eu{F}(S)\tir &\lim(\eu{F}(S\inv{\varpi})\times \eu{F}(\hat{S})\dutir &\eu{F}(\hat{S}\inv{\varpi})).\end{teqn*} 
The ring $S$, being flat over $R$, has bounded $\varpi$-torsion, consequently, so does $\hat{S}$ and, in fact, $\hat{S}\langle\varpi^{\infty}\rangle = S\langle\varpi^\infty\rangle$. This implies that Beauville--Laszlo gluing condition is satisfied ($S\to\hat{S}$ is a `gluable pair' as in \cite[Thm.~6.4]{bhatt-mathew-arc}), proving the required equivalence.
\par For the final assertion of the claim, the three functors in question are arc sheaves locally of finite presentation (\Cref{prop:fet-arc-descent}) and \Cref{lem:lifting-et} shows that, in fact, there is an isomorphism $\eu{F}(A)\cong\eu{F}(A/I)$ for any henselian pair $(A,I)$. 
\end{proof}
\bt\label{thm:tilting-U-fet}
Let \letperf, and let \lettilt. Suppose \bcpair. Then there are compatible equivalences, functorial in $A$ and $U$ and compatible with orthogonality relation on idempotents, 
\begin{gather}
	\begin{tikzpicture}
		\label{diag:tilting-U-idem}
			\node (F) at (-2,0) {\phantom{To ensure align}};
			\node (U) at (0,0) {$\idem(U)\cong\idem(\tilt{U}),$};
			\node (A) at (2,0) { \hspace{1.75cm}and};
	\end{tikzpicture}
	\\ 
	\begin{tikzcd}
		\label{diag:tilting-U-fet}
				\fet/U\cong\fet/\tilt{U}.
	\end{tikzcd}
\end{gather} 
Moreover, given a commutative, finite \'etale $U$-group scheme $G$, there are a commutative, finite \'etale $\tilt{U}$-group scheme $\tilt{G}$ obtained by \eqref{diag:tilting-U-fet} and a functorial in $A$, $U$ and $G$ equivalence \bd\label{diag:tilting-U-et}\R{\Gamma}_{\etale}(U, G)\cong\R{\Gamma}_{\etale}(\tilt{U}, \tilt{G}).\ed
\et
We note that the inspiration for defining the `tilt' $\tilt{U}$ of $U$ comes from the tilting functor in the adic theory of perfectoid spaces. In \ccite{ces-scholze}{Thm.~2.2.7}, the authors showed the equivalence \eqref{diag:tilting-U-idem} and \eqref{diag:tilting-U-et}, the latter in the case of constant coefficients, that is, by replacing both $G$ and $\tilt{G}$ by an abstract abelian group. The proof given below is an adaptation of the proof of loc.~cit. 
\begin{proof}[Proof of \Cref{thm:tilting-U-fet}]

For any scheme $X$, the commutative, finite \'etale $X$-group schemes are commutative group objects in the category of finite \'etale $X$-schemes. Thus, given a commutative, finite \'etale $U$-group scheme $G$, the identification \eqref{diag:tilting-U-fet}, functorially in $A$ and $U$, identifies $G$ with a commutative, finite \'etale $\tilt{U}$-group scheme.
\par Our strategy is to find a suitable open covering of $U$ which has a corresponding open covering of $\tilt{U}$: precisely, we would like each open $V$ in the covering of $U$ to have a corresponding open $\tilt{V}$ (i.e., such that $\spec A\inv{\varpi}\subset V$ and $\spec \tilt{A}\inv{\tilt{\varpi}}\subset\tilt{V}$, and $\spec A\setminus V$ identifies with $\spec\tilt{A}\setminus\tilt{V}$ under \eqref{diag:identification}) in the covering of $\tilt{U}$; and, therefore, by Zariski descent of idempotents (resp., of $\fet$, resp., of \'etale cohomology), we will reduce to producing the functorial equivalence \eqref{diag:tilting-U-idem} (resp., \eqref{diag:tilting-U-fet}, resp., \eqref{diag:tilting-U-et}) for each of the open $V$ in covering of $U$. In fact, it is possible to obtain each open $V$ such that it is covered by $\spec A[\frac{1}{\varpi}]$ and an open $\spec B\subset\spec A$, and such that $\tilt{V}$ is covered by $\spec \tilt{A}[\frac{1}{\tilt{\varpi}}]$ and an open $\spec B'\subset\spec \tilt{A}$ in the following way: covering $Z$ by principal open subsets $\textstyle\{\spec A[\frac{1}{f}]\subset U\}$, we can find elements $\tilt{f}\in\tilt{A}$ (guaranteed by \Cref{defn:integral-perfd}\eqref{pt:2} and \eqref{diag:lem-3.2}), for each $f$, such that $f^{\flat\sharp}\equiv f\pmod{\varpi}$, and such that the opens $\{\spec \tilt{A}[\frac{1}{\tilt{f}}]\subset\tilt{U}\}$ cover $\tilt{Z}$. Indeed, fixing one such cover of $Z$, for each open $\spec A\inv{f}$ in the cover, we may take $B=A\inv{f}$ and $B'=\tilt{A}\inv{\tilt{f}}$, and 
\begin{center}\begin{tikzpicture}
\node (V) at (0,0) {$\textstyle V=\spec A\inv{\varpi}\cup\spec B$};
\node (Vtilt) at (0.1,-1) {$\textstyle \tilt{V}=\spec\tilt{A}\inv{\tilt{\varpi}}\cup\spec B'.$};
\node (and) at (3,0.05) {and};
\end{tikzpicture}\end{center}
\par Henceforth, we asssume that $U$ is the union of $\spec A\inv{\varpi}$ and $\spec B$ (resp., $\tilt{U}$ is the union of $\spec\tilt{A}\inv{\tilt{\varpi}}$ and $\spec B'$) for some open $\spec B\subset\spec A$ (resp., $\spec B'\subset\spec\tilt{A}$). Our construction implies that $B/\varpi\cong B'/\tilt{\varpi}$, and therefore, the $\varpi$-adic completion $\hat{B}$ of $B$, which is perfectoid thanks to \cite[Cor.~2.1.6]{ces-scholze}, has tilt $\hat{B'}$, the $\tilt{\varpi}$-adic completion of $B'$ (which is also perfectoid thanks to loc.~cit.). 
The Beauville--Laszlo gluing (\Cref{prop:beauville-laszlo}) 
applies by using $R=A$ and $S=\varpi$-henselisation of $B$ (resp., by using $R=\tilt{A}$ and $S=\tilt{\varpi}$-henselisation of $B'$), reducing the task of producing a functorial equivalence \eqref{diag:tilting-U-idem} (resp., \eqref{diag:tilting-U-fet}, resp., \eqref{diag:tilting-U-et}) to the following two cases, namely, when $U=\spec A[\frac{1}{\varpi}]$ and $\tilt{U}=\spec \tilt{A}[\frac{1}{\tilt{\varpi}}]$ and when $U=\spec A$ and $\tilt{U}=\spec \tilt{A}$ (therefore, in particular, implying the case when $U=\spec \hat{B}$ and $\tilt{U}=\spec \hat{B'}$). The latter case is easier, because we can use \eqref{diag:identification} and \Cref{lem:lifting-et} (this is applicable since $A$ and $\tilt{A}$ are $\varpi$ and $\tilt{\varpi}$-henselian respectively). In the case of idempotents, the first case follows from \eqref{diag:lem-3.2}. 
\par To deal with the cases of $\fet$ and the \'etale cohomology, we will need to use the descent results from \textsection\ref{section:arc-descent}. Let \bud\label{diag:thm-tilting1} A\arrow[r]& A_0\dutir&A_1\arrow[r, shift left=0.2cm]\arrow[r]\arrow[r, shift right =0.2cm]&\cdots\eud be a $\varpi$-complete arc hypercover supplied by \Cref{lem:arc-hypercover2} and \bud\label{diag:thm-tilting2}\tilt{A}\arrow[r]& \tilt{A_0}\dutir&\tilt{A_1}\arrow[r, shift left=0.2cm]\arrow[r]\arrow[r, shift right =0.2cm]&\cdots.\eud be the tilt $\tilt{\varpi}$-complete arc hypercover.
By \Cref{rem:arc-w-descent}, the functor on $\varpi$-complete $A$-algebras taking such an $A$-algebra $\hjbrl{A}$ to the category of finite \'etale $\hjbrl{A}[\frac{1}{\varpi}]$-algebras (resp., to $\R{\Gamma}_{\etale}(\hjbrl{A}[\frac{1}{\varpi}],G)$) satisfies $\varpi$-complete arc hyperdescent and the functor on $\tilt{\varpi}$-complete $\tilt{A}$-algebras taking such an $\tilt{A}$-algebra $\hjbrl{A}$ to the category of finite \'etale $\hjbrl{A}[\frac{1}{\tilt{\varpi}}]$-algebras (resp., to $\R{\Gamma}_{\etale}(\hjbrl{A}[\frac{1}{\tilt{\varpi}}],\tilt{G})$) satisfies $\tilt{\varpi}$-complete arc hyperdescent. Hence, to show that there is a functorial equivalence $\fet/A[\frac{1}{\varpi}]\cong\fet/\tilt{A}[\frac{1}{\tilt{\varpi}}]$ (resp., $\R{\Gamma}_{\etale}(A[\frac{1}{\varpi}],G)\cong\R{\Gamma}_{\etale}(\tilt{A}[\frac{1}{\tilt{\varpi}}],\tilt{G})$), it is enough to exhibit functorial equivalences for all $i$, 
\begin{center}\begin{tikzpicture}
\node (fet) at (0,0) {$\fet/A_i[\frac{1}{\varpi}]\cong\fet/\tilt{A_i}[\frac{1}{\tilt{\varpi}}]$};
\node (et) at (-0.38,-1) {$\text{ (resp., }\R{\Gamma}_{\etale}(A_i[\frac{1}{\varpi}],G)\cong\R{\Gamma}_{\etale}(\tilt{A_i}[\frac{1}{\tilt{\varpi}}],\tilt{G})).$};
\end{tikzpicture}\end{center}
\par Because of the nature of the rings $A_i$, it is enough to establish functorial equivalences \eqref{diag:tilting-U-fet} and \eqref{diag:tilting-U-et} in the case when $\ma{R}=\prod_{i}V_i$, where $V_i$ are $\varpi$-adically complete \aic. 
By, for example \cite[Prop.~2.1.8]{ces-scholze}, $V_i$ are perfectoid and by \Cref{prop:alg-tilting}, $\tilt{V_i}$ are $\tilt{\varpi}$-adically complete valuation rings of rank $\le 1$ with algebraically closed fraction fields; consequently, $\ma{R}$ is perfectoid with tilt $\tilt{\ma{R}}$ due to op.~cit.~Prop. 2.1.11(d). 
By \Cref{lem:aic-v-fet}, the finite étale schemes over $R$ (resp., over $\tilt{R}$)
correspond to disjoint unions of subsets of $\spec R$ (resp., of $\spec\tilt{R}$) which are both open and closed, which, by \cite[\href{https://stacks.math.columbia.edu/tag/00EE}{Tag 00EE}]{stacks-project}, correspond to finite collections of idempotents of $R$ (resp., of $\tilt{R}$); whence, the functorial equivalence \eqref{diag:tilting-U-idem} implies \eqref{diag:tilting-U-fet}.
Again, by \cite[Lem. 2.2.9]{ces-scholze}, the ring $R$ (resp., the ring $\tilt{R}$) has no non-split étale covers, and therefore $\R{\Gamma}_{\etale}(R, G)$ (resp., $\R{\Gamma}_{\etale}(\tilt{R}, \tilt{G})$) is concentrated in degree $0$; thus, the cohomology $\R{\Gamma}_{\etale}(R,G)\cong H^0(R,G)$ (resp., the cohomology $\R{\Gamma}_{\etale}(\tilt{R},\tilt{G})\cong H^0(\tilt{R},\tilt{G})$) can be identified with the group of sections $R\to G$ (resp., sections $\tilt{R}\to\tilt{G}$). Consequently, the full faithfulness of the functorial equivalence \eqref{diag:tilting-U-fet} implies \eqref{diag:tilting-U-et}.
\end{proof}
In the above proof, we could have replaced the $\varpi$-complete arc hyperdescent with the `$\varpi$-complete $v$'-hyperdescent. Indeed, the proof of \ccite{ces-scholze}{Lem.~2.2.3} can be tailored to produce a `$\varpi$-complete $v$'-cover in \Cref{lem:arc-hypercover2}. However, it does not save us much work if we just prove the `$\varpi$ complete $v$'-hyperdescent for the functors in \Cref{prop:fet-arc-descent} and \Cref{rem:arc-w-descent}.  
\br Let $R$ be a perfectoid Banach $K$-algebra as in \cite[Defn.~5.1]{scholze-thesis}, where $K$ is a perfectoid field as in op.~cit.~Defn.~3.1, (resp., perfectoid Banach $\bb{Q}_p$-algebra as in \cite[Defn.~3.6.1]{kedlaya-liu}). 
Then, the ring $R^\circ$ of power-bounded elements is a perfectoid ring (cf.~\cite[Lem.~3.20]{bms1}). 
Choosing an element $\varpi\in R^\circ$ such that $\varpi^p\mid p$ and that $R^\circ$ is $\varpi$-adically complete, and choosing $\tilt{\varpi}\in R^{\flat\circ}$ such that $\varpi^{\flat\sharp}$ is a unit multiple of $\varpi$, we have, by \Cref{thm:tilting-U-fet}, an equivalence \bud\textstyle\fet/R^\circ[\frac{1}{\varpi}]\cong\fet/{R^{\flat\circ}}[\frac{1}{\tilt{\varpi}}].\eud 
\er
\bc[\ccite{ces-brauer}{Thm. 4.10}]\label{cor:main}
Let $A$ be a $\bb{Z}_p$-algebra such that it is a perfectoid ring with an element $\varpi\in A$ such that $\varpi^p\mid p$ and that $A$ is $\varpi$-adically complete. Then, for a commutative, finite étale $A[\frac{1}{\varpi}]$-group scheme $G$ of $p$-power order, we have, for all $i\ge 2$, \bd\label{diag:cor-main}\textstyle H^i_{\etale}(A[\frac{1}{\varpi}],G)=0.\ed In particular, for a commutative, finite étale $A[\frac{1}{p}]$-group scheme $G$ of $p$-power order, we have, for all $i\ge 2$, \bud\textstyle H^i_{\etale}(A[\frac{1}{p}],G)=0.\eud
\ec
\begin{proof}\let\qed\relax
The second vanishing follows from \eqref{diag:cor-main} because there exists a $\varpi\in A$ such that $\varpi^p$ is a unit multiple of $p$, and $A[\frac{1}{\varpi}]=A[\frac{1}{p}]$. Indeed, \cite[Lem.~3.9]{bms1} implies that there exists a unit $v\in A$ such that $vp$ admits compatible $p$-power roots, consequently, we can take $\varpi\in A$ such that $\varpi^p=vp$. 
\par Letting $\tilt{\varpi}\in\tilt{A}$ be such that $\varpi^{\flat\sharp}=vp$, \Cref{thm:tilting-U-fet} implies that we have an isomorphism \bud\textstyle\R{\Gamma}_{\etale}(A[\frac{1}{\varpi}],G)\cong\R{\Gamma}_{\etale}(\tilt{A}[\frac{1}{\tilt{\varpi}}],\tilt{G}).\eud It suffices to show the cohomology vanishing of the second complex, which reduces us to prove \eqref{diag:cor-main} for the perfect $\bb{F}_p$-algebra $\tilt{A}$. Due to a limit argument, \cite[Ex.~X Thm.~5.1]{sga4iii}, which shows that the $p$-cohomological degree of any noetherian $\bb{F}_p$-algebra is $\le 1$, shows that for all $i\ge 2$, \bd\tag*{$\qedsymbol$}\textstyle H^i_{\etale}(\tilt{A}[\frac{1}{\tilt{\varpi}}], \tilt{G})=0.\ed
\end{proof}

\printbibliography

@misc{bhatt-perfectoid,
	title		= {Lecture notes for a class on perfectoid spaces},
	author 		= {Bhatt, Bhargav},
	year		= {2017},
	url			= {http://www-personal.umich.edu/~bhattb/teaching/mat679w17/lectures.pdf},
}

@misc{bhatt-scholze,
	title		= {Prisms and Prismatic Cohomology},
	author		= {Bhargav Bhatt and Peter Scholze},
	year		= {2019},
	eprint		= {1905.08229},
	archivePrefix	= {arXiv},
	version={v3},
	primaryClass	= {math.AG}
}

@article {bhatt-mathew-arc,
    AUTHOR = {Bhatt, Bhargav and Mathew, Akhil},
     TITLE = {The arc-topology},
   JOURNAL = {Duke Math. J.},
  FJOURNAL = {Duke Mathematical Journal},
    VOLUME = {170},
      YEAR = {2021},
    NUMBER = {9},
     PAGES = {1899--1988},
      ISSN = {0012-7094},
   MRCLASS = {14F20 (14F06 14G22)},
  MRNUMBER = {4278670},
       DOI = {10.1215/00127094-2020-0088},
}

@misc{bouthier-ces,
      title={Torsors on loop groups and the Hitchin fibration}, 
      author={Alexis Bouthier and Česnavičius, Kęstutis},
      year={2020},
      version={v3},
      eprint={1908.07480},
      archivePrefix={arXiv},
      primaryClass={math.AG}
}

@article{bms1,
	AUTHOR = {Bhatt, Bhargav and Morrow, Matthew and Scholze, Peter},
     TITLE = {Integral {$p$}-adic {H}odge theory},
   JOURNAL = {Publ. Math. Inst. Hautes \'{E}tudes Sci.},
  FJOURNAL = {Publications Math\'{e}matiques. Institut de Hautes \'{E}tudes
              Scientifiques},
    VOLUME = {128},
      YEAR = {2018},
     PAGES = {219--397},
      ISSN = {0073-8301},
   MRCLASS = {14F30},
  MRNUMBER = {3905467},
MRREVIEWER = {Daniel Robert Gulotta},
       DOI = {10.1007/s10240-019-00102-z},
}

@article{ces-brauer,
    AUTHOR = {Česnavičius, Kęstutis},
     TITLE = {Purity for the {B}rauer group},
   JOURNAL = {Duke Math. J.},
  FJOURNAL = {Duke Mathematical Journal},
    VOLUME = {168},
      YEAR = {2019},
    NUMBER = {8},
     PAGES = {1461--1486},
      ISSN = {0012-7094},
   MRCLASS = {14F22 (14F20 14G22 16K50)},
  MRNUMBER = {3959863},
MRREVIEWER = {J\"{o}rg Jahnel},
       DOI = {10.1215/00127094-2018-0057},
}

@misc{ces-scholze,
      title={Purity for flat cohomology}, 
      author={Kęstutis Česnavičius and Peter Scholze},
      year={2021},
      eprint={1912.10932},
      archivePrefix={arXiv},
      version={v2},
      primaryClass={math.AG}
}

@article {ehik,
    AUTHOR = {Elmanto, Elden and Hoyois, Marc and Iwasa, Ryomei and Kelly,
              Shane},
     TITLE = {Cdh descent, cdarc descent, and {M}ilnor excision},
   JOURNAL = {Math. Ann.},
  FJOURNAL = {Mathematische Annalen},
    VOLUME = {379},
      YEAR = {2021},
    NUMBER = {3-4},
     PAGES = {1011--1045},
      ISSN = {0025-5831},
   MRCLASS = {14F42 (14A20)},
  MRNUMBER = {4238259},
       DOI = {10.1007/s00208-020-02083-5},
}

@article{gabber_1994, 
	AUTHOR = {Gabber, Ofer},
     TITLE = {Affine analog of the proper base change theorem},
   JOURNAL = {Israel J. Math.},
  FJOURNAL = {Israel Journal of Mathematics},
    VOLUME = {87},
      YEAR = {1994},
    NUMBER = {1-3},
     PAGES = {325--335},
      ISSN = {0021-2172},
   MRCLASS = {14F20 (13J15)},
  MRNUMBER = {1286833},
MRREVIEWER = {Gennady Lyubeznik},
       DOI = {10.1007/BF02773001},
}

@book {huber96,
    AUTHOR = {Huber, Roland},
     TITLE = {\'{E}tale cohomology of rigid analytic varieties and adic spaces},
    SERIES = {Aspects of Mathematics, E30},
 PUBLISHER = {Friedr. Vieweg \& Sohn, Braunschweig},
      YEAR = {1996},
     PAGES = {x+450},
      ISBN = {3-528-06794-2},
   MRCLASS = {14G22 (14F20)},
  MRNUMBER = {1734903},
MRREVIEWER = {Lorenzo Ramero},
       DOI = {10.1007/978-3-663-09991-8},
}

@article {kedlaya-liu,
    AUTHOR = {Kedlaya, Kiran S. and Liu, Ruochuan},
     TITLE = {Relative {$p$}-adic {H}odge theory: foundations},
   JOURNAL = {Ast\'{e}risque},
  FJOURNAL = {Ast\'{e}risque},
    NUMBER = {371},
      YEAR = {2015},
     PAGES = {239},
      ISSN = {0303-1179},
      ISBN = {978-2-85629-807-7},
   MRCLASS = {14C30 (14G22)},
  MRNUMBER = {3379653},
MRREVIEWER = {Giovanni Rosso},
}

@unpublished{master-thesis,
	title={Tiltings Correspondences on Perfectoid Rings},
	author={Kundu, Arnab},
	year={2020},
	type={Masters' Thesis},
}

@unpublished{lurie-algebra,
	title={Higher Algebra},
	author={Lurie, Jacob},
	year={2017},
	type={book},
	shorthand = {HA},
	url={https://www.math.ias.edu/~lurie/papers/HA.pdf},
}

@book {lurie-htt,
 shorthand = {HTT},
    AUTHOR = {Lurie, Jacob},
     TITLE = {Higher topos theory},
    SERIES = {Annals of Mathematics Studies},
    VOLUME = {170},
 PUBLISHER = {Princeton University Press, Princeton, NJ},
      YEAR = {2009},
     PAGES = {xviii+925},
      ISBN = {978-0-691-14049-0; 0-691-14049-9},
   MRCLASS = {18-02 (18B25 18E35 18G30 18G55 55U40)},
  MRNUMBER = {2522659},
MRREVIEWER = {Mark Hovey},
       DOI = {10.1515/9781400830558},
}

@misc{mathew,
      title={Faithfully flat descent of almost perfect complexes in rigid geometry}, 
      author={Akhil Mathew},
      year={2020},
      eprint={1912.10968},
      archivePrefix={arXiv},
      version={v2},
      primaryClass={math.AG}
}

@article {rydh,
    AUTHOR = {Rydh, David},
     TITLE = {Submersions and effective descent of \'{e}tale morphisms},
   JOURNAL = {Bull. Soc. Math. France},
  FJOURNAL = {Bulletin de la Soci\'{e}t\'{e} Math\'{e}matique de France},
    VOLUME = {138},
      YEAR = {2010},
    NUMBER = {2},
     PAGES = {181--230},
      ISSN = {0037-9484},
   MRCLASS = {14A15 (13B22 13B40 14F20 14F43)},
  MRNUMBER = {2679038},
MRREVIEWER = {Liam O'Carroll},
       DOI = {10.24033/bsmf.2588},
}

@article{scholze-thesis,
	AUTHOR = {Scholze, Peter},
     TITLE = {Perfectoid spaces},
   JOURNAL = {Publ. Math. Inst. Hautes \'{E}tudes Sci.},
  FJOURNAL = {Publications Math\'{e}matiques. Institut de Hautes \'{E}tudes
              Scientifiques},
    VOLUME = {116},
      YEAR = {2012},
     PAGES = {245--313},
      ISSN = {0073-8301},
   MRCLASS = {14G99},
  MRNUMBER = {3090258},
MRREVIEWER = {Jean-Marc Fontaine},
       DOI = {10.1007/s10240-012-0042-x},
}

@book {sga1,
     TITLE = {Rev\^{e}tements \'{e}tales et groupe fondamental ({SGA} 1)},
    SERIES = {Documents Math\'{e}matiques (Paris) [Mathematical Documents
              (Paris)]},
    VOLUME = {3},
    shorthand = {SGA1},
      NOTE = {S\'{e}minaire de g\'{e}om\'{e}trie alg\'{e}brique du Bois Marie 1960--61.
              [Algebraic Geometry Seminar of Bois Marie 1960-61],
              Directed by A. Grothendieck,
              With two papers by M. Raynaud,
              Updated and annotated reprint of the 1971 original [Lecture
              Notes in Math., 224, Springer, Berlin;  MR0354651 (50
              \#7129)]},
 PUBLISHER = {Soci\'{e}t\'{e} Math\'{e}matique de France, Paris},
      YEAR = {2003},
     PAGES = {xviii+327},
      ISBN = {2-85629-141-4},
   MRCLASS = {14E20 (14-06 14F35)},
  MRNUMBER = {2017446},
}

@book {sga4iii,
     TITLE = {Th\'{e}orie des topos et cohomologie \'{e}tale des sch\'{e}mas. {T}ome 3},
     shorthand = {SGA4 III},
    SERIES = {Lecture Notes in Mathematics, Vol. 305},
      NOTE = {S\'{e}minaire de G\'{e}om\'{e}trie Alg\'{e}brique du Bois-Marie 1963--1964
              (SGA 4),
              Dirig\'{e} par M. Artin, A. Grothendieck et J. L. Verdier. Avec la
              collaboration de P. Deligne et B. Saint-Donat},
 PUBLISHER = {Springer-Verlag, Berlin-New York},
      YEAR = {1973},
     PAGES = {vi+640},
   MRCLASS = {14-06},
  MRNUMBER = {0354654},
}

@misc{stacks-project,
	author       	= {The {Stacks project authors}},
	title        	= {The Stacks project},
	howpublished 	= {\url{https://stacks.math.columbia.edu}},
	year         	= {2020},
}
\end{document}